\documentclass[a4paper,12pt]{amsart}
\usepackage{xr-hyper}
\usepackage{amssymb,amscd,amsmath,a4wide,graphicx,stmaryrd,fullpage,setspace,microtype,textcomp,tikz,enumitem,accents,mathtools,abraces,cite}
\usepackage{hyperref}
\usepackage{float}
\usepackage[all,cmtip]{xy}
\title{Erratum to ``Homotopy theory of Moore flows I''}

\author[P. Gaucher]{Philippe Gaucher}

\address{Universit\'e Paris Cit\'e, CNRS, IRIF, F-75013, Paris, France}

\urladdr{http://www.irif.fr/{\~{}}gaucher} 

\makeatletter
\@namedef{subjclassname@2020}{\textup{2020} Mathematics Subject Classification}
\makeatother


\swapnumbers

\newcommand{\K}{\mathcal{K}}

\newcommand{\p}{\times}

\newtheorem*{thmN}{Theorem}
\newtheorem{thm}{Theorem}[section]
\newtheorem{prop}[thm]{Proposition}
\newtheorem{lem}[thm]{Lemma}

\newtheorem{cor}[thm]{Corollary}

\newtheorem{defn}[thm]{Definition}
\newtheorem{nota}[thm]{Notation}

\newtheorem{defnot}[thm]{Definition and notation}
\newcommand{\bd}{\begin{defn}}
\newcommand{\ed}{\end{defn}}
\newcommand{\bdn}{\begin{defnot}}
\newcommand{\edn}{\end{defnot}}
\newcommand{\bp}{\begin{prop}}
\newcommand{\ep}{\end{prop}}
\newcommand{\bth}{\begin{thm}}
\renewcommand{\eth}{\end{thm}}
\newcommand{\bpf}{\begin{proof}}
\newcommand{\epf}{\end{proof}}
\newcommand{\bc}{\begin{cor}}
\newcommand{\ec}{\end{cor}}

\newcommand{\fL}[1]{\ar@{->}[ll]_-{#1}}
\newcommand{\fR}[1]{\ar@{->}[rr]^-{#1}}
\newcommand{\fRr}[1]{\ar@{->}[rrr]^-{#1}}
\newcommand{\fD}[1]{\ar@{->}[dd]_-{#1}}
\newcommand{\fU}[1]{\ar@{->}[uu]^-{#1}}
\newcommand{\f}[2]{\ar@{->}[#1]|{#2}}
\newcommand{\ff}[2]{\ar@2{->}[#1]|{#2}}
\newcommand{\frr}[1]{\ar@{->}[rrrr]^-{#1}}

\newcommand{\fl}[1]{\ar@{->}[l]_-{#1}}
\newcommand{\fr}[1]{\ar@{->}[r]^-{#1}}
\newcommand{\fd}[1]{\ar@{->}[d]_-{#1}}
\newcommand{\fu}[1]{\ar@{->}[u]^-{#1}}

\renewcommand{\top}{{\mathbf{Top}}}

\newcommand{\iso}{\cong}

\renewcommand{\leq}{\leqslant}
\renewcommand{\geq}{\geqslant}

\newcommand{\topdgr}{[\mathcal{P}^{op},\top]}

\def\cartesien{%
  \ar@{-}[]+R+<6pt,-2pt>;[]+RD+<6pt,-6pt>%
  \ar@{-}[]+D+<2pt,-6pt>;[]+RD+<6pt,-6pt>%
}
\def\cocartesien{%
  \ar@{-}[]+L+<-6pt,+2pt>;[]+LU+<-6pt,+6pt>%
  \ar@{-}[]+U+<-2pt,+6pt>;[]+LU+<-6pt,+6pt>%
}
\def\hocartesien{%
  \ar@{-}[]+R+<6pt,-2pt>;[]+RD+<6pt,-6pt>_{h}%
  \ar@{-}[]+D+<2pt,-6pt>;[]+RD+<6pt,-6pt>%
}
\def\hococartesien{%
  \ar@{-}[]+L+<-6pt,+2pt>;[]+LU+<-6pt,+6pt>_{h}%
  \ar@{-}[]+U+<-2pt,+6pt>;[]+LU+<-6pt,+6pt>%
}

\newcommand{\brm}[1]{\rm{\mathbf{#1}}}

\newcommand{\ttop}{{\brm{TOP}}}

\DeclareMathOperator{\id}{Id}

\DeclareMathOperator{\Obj}{Obj}

\newcommand{\liminj}{\varinjlim}

\makeatletter
\def\varholim@#1#2{%
  \vtop{\m@th\ialign{##\cr
    \hfil$#1\operator@font holim$\hfil\cr
    \noalign{\nointerlineskip\kern1.5\ex@}#2\cr
    \noalign{\nointerlineskip\kern-\ex@}\cr}}%
}
\def\holimproj{%
  \mathop{\mathpalette\varholim@{\leftarrowfill@\textstyle}}\nmlimits@
}
\def\holiminj{%
  \mathop{\mathpalette\varholim@{\rightarrowfill@\textstyle}}\nmlimits@
}
\makeatother

\makeatletter
\newcommand*{\@opargbegintheorem}[3]{\trivlist
	\item[\hskip \labelsep{\bfseries #1\ #2}] \textbf{(#3)}\ \itshape}
\makeatother
\date{March 14, 2022}

\newcommand{\adj}[4]{\xymatrix@1{{#1}\ar@/^0.8em/[r]^-{#2} \ar@{}[r]|-{\perp} & \ar@/^0.8em/[l]^-{#3} {#4}}}

\newcommand{\ot}{\otimes}

\begin{document}

\begin{abstract} 
	The notion of reparametrization category is incorrectly axiomatized and it must be adjusted. It is proved that for a general reparametrization category $\mathcal{P}$, the tensor product of $\mathcal{P}$-spaces yields a biclosed semimonoidal structure. It is also described some kind of objectwise braiding for $\mathcal{G}$-spaces.
\end{abstract}

\maketitle

\tableofcontents

\section{Introduction}

\subsection*{Presentation} The notion of reparametrization category introduced in \cite{Moore1} is incorrectly axiomatized. The reparametrization categories $(\mathcal{G},+)$ and $(\mathcal{M},+)$ are not symmetric indeed. Moreover, the third axiom of reparametrization category is slightly modified to obtain the expected result for the tensor product of two constant $\mathcal{P}$-spaces in full generality. It also enables us to write a short proof of the pentagon axiom. The main theorem is:

\begin{thmN} (Proposition~\ref{asso} and Theorem~\ref{closedsemi})
	For any reparametrization category $\mathcal{P}$, the tensor product of $\mathcal{P}$-spaces yields a biclosed  semimonoidal structure.
\end{thmN}

The semimonoidal category of $\mathcal{G}$-spaces still has some kind of objectwise braiding which is formalized in Theorem~\ref{non-natural-homeo}. This fact is specific to $\mathcal{G}$-spaces. It is used nowhere in \cite{Moore1,Moore2}.

\begin{thmN} (Theorem~\ref{non-natural-homeo})
There is a homeomorphism \[B: (D\ot E)(L) \longrightarrow (E\ot D)(L)\] for all $L>0$ and all $\mathcal{G}$-spaces $D$ and $E$ which is not natural with respect to $L>0$.
\end{thmN}

\subsection*{Outline of the note} In Section~\ref{reparam}, the notion of reparametrization category is adjusted. In Section~\ref{rep}, the corrections are listed. The absence of braiding forces us to relocate some parameters $\ell$ in the calculations, and also to replace the shift operator $s_\ell$ either by the \textit{left shift} $s^L_\ell$ (see Proposition~\ref{dec_left}) or by the \textit{right shift} $s^R_\ell$ (see Proposition~\ref{dec_right}). Finally, Section~\ref{non-natural} gives an explicit description of a homeomorphism $(D\ot E)(L) \iso (E\ot D)(L)$ for all $L>0$ and for all $\mathcal{G}$-spaces $D$ and $E$ which is not natural with respect to $L>0$. 

\subsection*{Prerequisites and notations} We refer to \cite{Moore1} for the notations and for the full categorical argumentations. We refer to \cite{Moore2} for the full topological argumentations.

\section{Adjustment}
\label{reparam}

\bd A {\rm semimonoidal category} $(\K,\ot)$ is a category $\K$ equipped with a functor $\ot:\K\p \K\to \K$ together with a natural isomorphism $a_{x,y,z}: (x\ot y) \ot z \to x \ot (y\ot z)$ called the {\rm associator} satisfying the pentagon axiom.
\ed

\bd A semimonoidal category $(\K,\ot)$ is {\rm enriched} (all enriched categories are enriched over $\top$) if the category $\K$ is enriched and if the set map \[\K(a,b)\p \K(c,d) \longrightarrow \K(a\ot c,b\ot d)\] is continuous for all objects $a,b,c,d\in \Obj(\K)$. 
\ed

\bd \label{def-reparam}
A {\rm reparametrization category} $(\mathcal{P},\ot)$ is a small enriched semimonoidal category satisfying the following additional properties: 
\begin{enumerate}
	\item The semimonoidal structure is strict, i.e. the associator is the identity.
	\item All spaces of maps $\mathcal{P}(\ell,\ell')$ for all objects $\ell$ and $\ell'$ of $\mathcal{P}$ are contractible. 
	\item For all maps $\phi:\ell\to \ell'$ of $\mathcal{P}$, for all $\ell'_1,\ell'_2\in \Obj(\mathcal{P})$ such that $\ell'_1\ot\ell'_2=\ell'$, there exist two maps $\phi_1:\ell_1\to \ell'_1$ and $\phi_2:\ell_2\to \ell'_2$ of $\mathcal{P}$ such that $\phi=\phi_1 \ot \phi_2 : \ell_1\ot\ell_2 \to \ell'_1 \ot\ell'_2$ (which implies that $\ell_1 \ot \ell_2=\ell$). 
%	\item One has $\id_{\ell\ot \ell'}\ot \id_{\ell''}=\id_{\ell} \ot \id_{\ell' \ot \ell''}$ for all $\ell,\ell',\ell''\in\Obj(\mathcal{P})$.
\end{enumerate}
\ed 

\begin{nota}
	The notations $\ell,\ell',\ell_i,L,\dots$ mean objects of a reparametrization category $\mathcal{P}$. 
\end{nota}

\begin{nota}
	To stick to the intuition, we set $\ell+\ell' := \ell \ot \ell'$ for all $\ell,\ell'\in \Obj(\mathcal{P})$. Indeed, morally speaking, $\ell$ is the length of a path. 
\end{nota}

The enriched categories $(\mathcal{G},+)$ (Proposition~\ref{paramG}), $(\mathcal{M},+)$ \cite[Proposition~4.11]{Moore1} as well as the terminal category are examples of reparametrization categories. In the cases of $(\mathcal{G},+)$ and $(\mathcal{M},+)$, the functors $(\ell,\ell')\mapsto \ell+\ell'$ and $(\ell,\ell') \mapsto \ell'+\ell$ coincide on objects, but not on morphisms. The terminal category is a symmetric reparametrization category. We do not know if there exist symmetric reparametrization categories not equivalent to the terminal category. \cite[Proposition~5.8]{Moore1} must be replaced by the two propositions:

\bp \label{dec_left} (The left shift functor)
The following data assemble to an enriched functor $s^L_\ell:\mathcal{P}\to \mathcal{P}$: 
\[
\begin{cases}
s^L_\ell(\ell')= \ell + \ell'\\
s^L_\ell(\phi) = \id_\ell \ot \phi & \hbox{for a map $\phi:\ell'\to \ell''$.}
\end{cases}
\]
\ep

\bp \label{dec_right} (The right shift functor)
The following data assemble to an enriched functor $s^R_\ell:\mathcal{P}\to \mathcal{P}$: 
\[
\begin{cases}
s^R_\ell(\ell')= \ell' + \ell\\
s^R_\ell(\phi) =  \phi \ot \id_\ell & \hbox{for a map $\phi:\ell'\to \ell''$.}
\end{cases}
\]
\ep

For the convenience of the reader, we recall the 

\bd \label{pspace} \cite[Definition~5.1]{Moore1} An object of $\topdgr_0$ is called a {\rm $\mathcal{P}$-space}. Let $D$ be a $\mathcal{P}$-space. Let $\phi:\ell\to \ell'$ be a map of $\mathcal{P}$. Let $x\in D(\ell')$. We will use the notation 
\[
x.\phi := D(\phi)(x).
\]
\ed

\begin{nota}
	The two enriched functors $(s^L_\ell)^*$ and $(s^R_\ell)^*$ take a $\mathcal{P}$-space $D$ to $Ds^L_\ell$ and $Ds^R_\ell$ respectively.
\end{nota}

\section{Corrections}
\label{rep}

\begin{lem} \label{petitcalcul} (First replacement for \cite[Lemma~5.10]{Moore1})
	For all $\ell',\ell''\in \Obj(\mathcal{P})$, there is the isomorphism of $\mathcal{P}$-spaces (natural with respect to $\ell'$ and $\ell''$)
	\[
	\int^{\ell} \mathcal{P}(-,\ell+\ell') \p \mathcal{P}(\ell,\ell'') \iso \mathcal{P}(-,\ell''+\ell').
	\]
	The isomorphism takes the equivalence class of $(\psi,\phi)\in \mathcal{P}(-,\ell+\ell') \p \mathcal{P}(\ell,\ell'')$ to $(s^R_{\ell'})^*(\phi)\psi = (\phi\ot \id_{\ell'})\psi$. 
\end{lem}

\bpf
Pick a $\mathcal{P}$-space $D$. Then there is the sequence of homeomorphisms
\[\begin{aligned}
\topdgr\bigg(\int^{\ell} \mathcal{P}(-,\ell+\ell') \p \mathcal{P}(\ell,\ell''),D\bigg)
&\iso \int_{\ell} \topdgr\big(\mathcal{P}(-,\ell+\ell') \p \mathcal{P}(\ell,\ell''),D\big)\\
&\iso \int_{\ell} \ttop(\mathcal{P}(\ell,\ell''),D(\ell+\ell'))\\
&\iso \topdgr(\mathcal{P}(-,\ell''),(s^R_{\ell'})^*D)\\
&\iso D(\ell''+\ell')\\
&\iso \topdgr(\mathcal{P}(-,\ell''+\ell'),D).
\end{aligned}\]
The proof is complete thanks to the Yoneda lemma.
\epf

There is the following variation of Lemma~\ref{petitcalcul} which is also used below: 

\begin{lem} \label{petitcalcul2} (Second replacement for \cite[Lemma~5.10]{Moore1})
	For all $\ell',\ell''\in \Obj(\mathcal{P})$, there is the isomorphism of $\mathcal{P}$-spaces (natural with respect to $\ell'$ and $\ell''$)
	\[
	\int^{\ell} \mathcal{P}(-,\ell'+\ell) \p \mathcal{P}(\ell,\ell'') \iso \mathcal{P}(-,\ell'+\ell'').
	\]
	The isomorphism takes the equivalence class of $(\psi,\phi) \in \mathcal{P}(-,\ell'+\ell) \p \mathcal{P}(\ell,\ell'')$ to $(s^L_{\ell'})^*(\phi)\psi =(\id_{\ell'}\ot \phi)\psi$.
\end{lem}

\bpf
Pick a $\mathcal{P}$-space $D$. Then there is the sequence of homeomorphisms
\[\begin{aligned}
\topdgr\bigg(\int^{\ell} \mathcal{P}(-,\ell'+\ell) \p \mathcal{P}(\ell,\ell''),D\bigg)
&\iso \int_{\ell} \topdgr\big(\mathcal{P}(-,\ell'+\ell) \p \mathcal{P}(\ell,\ell''),D\big)\\
&\iso \int_{\ell} \ttop(\mathcal{P}(\ell,\ell''),D(\ell'+\ell))\\
&\iso \topdgr(\mathcal{P}(-,\ell''),(s^L_{\ell'})^*D)\\
&\iso D(\ell'+\ell'')\\
&\iso \topdgr(\mathcal{P}(-,\ell'+\ell''),D).
\end{aligned}\]
The proof is complete thanks to the Yoneda lemma.
\epf

\bp \label{underlying-set-tensor}
Let $D_1$ and $D_2$ be two $\mathcal{P}$-spaces and $L\in\Obj(\mathcal{P})$. Then the mapping $(x,y) \mapsto (\id,x,y)$ yields a surjective continuous map 
\[
\displaystyle\bigsqcup\limits_{\substack{(\ell_1,\ell_2)\\\ell_1+\ell_2=L}}  D_1(\ell_1)\p  D_2(\ell_2)\longrightarrow (D_1\ot D_2)(L).
\]
\ep

\bpf
Let $(\psi,x_1,x_2) \in \mathcal{P}(L,\ell_1+\ell_2) \p D_1(\ell_1)\p  D_2(\ell_2)$ be a representative of an element of $(D_1\ot D_2)(L)$. Then there exist two maps $\psi_i:\ell'_i\to \ell_i$ for $i=1,2$ such that $\psi=\psi_1\ot \psi_2$. By \cite[Corollary~5.13]{Moore1}, one has $(\psi,x_1,x_2)\sim (\id_L,x_1\psi_1,x_2\psi_2)$ in $(D_1\ot D_2)(L)$ and the proof is complete. 
\epf

\bp \label{asso} (Replacement for \cite[Proposition~5.11]{Moore1})
Let $D$ and $E$ be two $\mathcal{P}$-spaces. Let 
\[
D \ot E = \int^{(\ell_1,\ell_2)} \mathcal{P}(-,\ell_1+\ell_2) \p D(\ell_1) \p E(\ell_2).
\]
The pair $(\topdgr_0,\ot)$ is a semimonoidal category. 
%The functor \[\ot : \topdgr_0\p \topdgr_0 \to \topdgr_0\] induces a semimonoidal structure on $\topdgr_0$. 
\ep

\bpf Let $D_1,D_2,D_3$ be three $\mathcal{P}$-spaces. Let $a_{D_1,D_2,D_3}:(D_1\ot D_2) \ot D_3 \to D_1 \ot (D_2\ot D_3)$ be the composite of the isomorphisms (by using Lemma~\ref{petitcalcul} and Lemma~\ref{petitcalcul2})
\[\begin{aligned}
(D_1\ot &D_2) \ot D_3\\ &\iso  \int^{(\ell_1,\ell_2,\ell_3)} \bigg(\int^{\ell} \mathcal{P}(-,\ell+\ell_3) \p \mathcal{P}(\ell,\ell_1+\ell_2)\bigg)\p D_1(\ell_1) \p D_2(\ell_2) \p D_3(\ell_3)\\
&\iso \int^{(\ell_1,\ell_2,\ell_3)} \mathcal{P}(-,\ell_1+\ell_2+\ell_3) \p D_1(\ell_1) \p D_2(\ell_2) \p D_3(\ell_3)\\
&\iso \int^{(\ell_1,\ell_2,\ell_3)} \bigg(\int^{\ell} \mathcal{P}(-,\ell_1+\ell) \p \mathcal{P}(\ell,\ell_2+\ell_3)\bigg)\p D_1(\ell_1) \p D_2(\ell_2) \p D_3(\ell_3)\\
&\iso D_1 \ot (D_2\ot D_3).
\end{aligned}\]
Let $(\psi,(\phi,x_1,x_2),x_3)\in ((D\ot E)\ot F)(L)$ with $x_i\in D_i(\ell_i)$ for $i=1,2,3$ and $L\in\Obj(\mathcal{P})$. Write $\phi=\phi_1\ot\phi_2$ with $\phi_i:\ell'_i\to \ell_i$ for $i=1,2$ and $\psi=\psi_1\ot \psi_2 \ot \psi_3$ with $\psi_i:\ell''_i\to \ell'_i$ for $i=1,2,3$ with $\ell'_3=\ell_3$. In particular, $L=\ell''_1+\ell''_2+\ell''_3$. We obtain $(\psi,(\phi,x_1,x_2),x_3) \sim (\id_{L},(\id_{\ell''_1+\ell''_2},x_1\phi_1\psi_1,x_2\phi_2\psi_2),x_3\psi_3)$ in $((D\ot E)\ot F)(L)$. The above sequence of isomorphisms takes the equivalence class of $(\psi,(\phi,x_1,x_2),x_3)$ at first to the equivalence class of $((\id_{\ell''_1+\ell''_2}\ot \id_{\ell''_3})\id_{L},x_1\phi_1\psi_1,x_2\phi_2\psi_2,x_3\psi_3)$ by Lemma~\ref{petitcalcul}, and, since $(\id_{\ell''_1+\ell''_2}\ot \id_{\ell''_3})\id_{L} = (\id_{\ell''_1}\ot \id_{\ell''_2 +\ell''_3})\id_{L}$ and by Lemma~\ref{petitcalcul2}, to the equivalence class of $(\id_{L},x_1\phi_1\psi_1,(\id_{\ell''_2+\ell''_3},x_2\phi_2\psi_2,x_3\psi_3))$. We deduce that the associator $a_{D,E,F}:(D\ot E) \ot F \to D \ot (E\ot F)$ satisfies the pentagon axiom using Proposition~\ref{underlying-set-tensor}. 
\epf

\bth \label{closedsemi} (Replacement for \cite[Theorem~5.14]{Moore1})
Let $D$, $E$ and $F$ be three $\mathcal{P}$-spaces. Let  
\begin{align*}
&\{E,F\}_L:= \ell\mapsto \topdgr(E,(s^L_\ell)^*F),\\
&\{E,F\}_R:= \ell\mapsto \topdgr(E,(s^R_\ell)^*F).
\end{align*}
These yield two $\mathcal{P}$-spaces and there are the natural homeomorphisms
\begin{align*}
&\topdgr(D,\{E,F\}_L) \iso \topdgr(D\ot E,F),\\
&\topdgr(E,\{D,F\}_R) \iso \topdgr(D\ot E,F).
\end{align*}
Consequently, the functor \[\ot : \topdgr_0\p \topdgr_0 \to \topdgr_0\] induces a structure of biclosed semimonoidal structure on $\topdgr_0$.
\eth

\bpf
There are the sequences of natural homeomorphisms 
\[\begin{aligned}
\topdgr(D,\{E,F\}_L) & \iso \int_\ell \ttop\big(D(\ell),\topdgr(E,(s^L_\ell)^*F)\big)\\
&\iso \int_{(\ell,\ell')} \ttop\big(D(\ell),\ttop(E(\ell'),F(\ell+\ell'))\big) \\
&\iso \int_{(\ell,\ell')} \ttop\big(D(\ell)\p E(\ell'),F(\ell+\ell')\big) \\
&\iso \int_{(\ell,\ell')} \topdgr\big(\mathcal{P}(-,\ell+\ell') \p D(\ell)\p E(\ell'),F\big)\\
&\iso \topdgr(D\ot E,F)
\end{aligned}\]
and 
\[\begin{aligned}
\topdgr(E,\{D,F\}_R) & \iso \int_{\ell'} \ttop\big(E(\ell'),\topdgr(D,(s^R_{\ell'})^*F)\big)\\
&\iso \int_{(\ell,\ell')} \ttop\big(E(\ell'),\ttop(D(\ell),F(\ell+\ell'))\big) \\
&\iso \int_{(\ell,\ell')} \ttop\big(D(\ell)\p E(\ell'),F(\ell+\ell')\big) \\
&\iso \int_{(\ell,\ell')} \topdgr\big(\mathcal{P}(-,\ell+\ell') \p D(\ell)\p E(\ell'),F\big)\\
&\iso \topdgr(D\ot E,F).
\end{aligned}\]
\epf

\begin{nota}
	Let \[\mathbb{F}^{\mathcal{P}^{op}}_{\ell}U=\mathcal{P}(-,\ell)\p U \in \topdgr_0\] where $U$ is a topological space and where $\ell$ is an object of $\mathcal{P}$.
\end{nota}

\bp \label{Ftenseur} (Replacement for \cite[Proposition~5.16]{Moore1})
Let $U,U'$ be two topological spaces. Let $\ell,\ell'\in \Obj(\mathcal{P})$. There is the natural isomorphism of $\mathcal{P}$-spaces 
\[
\mathbb{F}^{\mathcal{P}^{op}}_{\ell}U \ot \mathbb{F}^{\mathcal{P}^{op}}_{\ell'}U' \iso \mathbb{F}^{\mathcal{P}^{op}}_{\ell+\ell'}(U\p U').
\]
\ep

\bpf
One has 
\[
\mathbb{F}^{\mathcal{P}^{op}}_{\ell}U \ot \mathbb{F}^{\mathcal{P}^{op}}_{\ell'}U' = \int^{(\ell_1,\ell_2)} \mathcal{P}(-,\ell_1+\ell_2) \p \mathcal{P}(\ell_1,\ell) \p \mathcal{P}(\ell_2,\ell') \p U \p U'.
\]
Using Lemma~\ref{petitcalcul2}, we obtain 
\[
\mathbb{F}^{\mathcal{P}^{op}}_{\ell}U \ot \mathbb{F}^{\mathcal{P}^{op}}_{\ell'}U' = \int^{\ell_1} \mathcal{P}(\ell_1,\ell) \p \mathcal{P}(-,\ell_1+\ell') \p U \p U'.
\]
Using Lemma~\ref{petitcalcul}, we obtain
\[
\mathbb{F}^{\mathcal{P}^{op}}_{\ell}U \ot \mathbb{F}^{\mathcal{P}^{op}}_{\ell'}U' = \mathcal{P}(-,\ell+\ell') \p U \p U'.
\]
\epf

\begin{nota}
	Let $U$ be a topological space. The constant $\mathcal{P}$-space $U$ is denoted by $\Delta_{\mathcal{P}^{op}}U$. 
\end{nota}

\bp \label{Ptenseur} (Replacement for \cite[Proposition~5.17]{Moore1})
Let $U$ and $U'$ be two topological spaces.  There is the natural isomorphism of $\mathcal{P}$-spaces
\[
\Delta_{\mathcal{P}^{op}} U \ot \Delta_{\mathcal{P}^{op}} U' \iso \Delta_{\mathcal{P}^{op}} (U \p U').
\]
\ep

\bpf
Since $\top$ is cartesian closed, it suffices to consider the case where $U=U'$ is a singleton. In that case, the topological space $(\Delta_{\mathcal{P}^{op}} U \ot \Delta_{\mathcal{P}^{op}} U')(L)$ is the quotient of the space \[\bigsqcup_{(\ell,\ell')} \mathcal{P}(L,\ell+\ell')\] by the identifications $(\phi_1 \ot \phi_2).\phi \sim \phi$. Let $\psi\in \mathcal{P}(L,\ell+\ell')$ for some $\ell,\ell'\in \Obj(\mathcal{P})$. By definition of a reparametrization category, write $\psi=\psi_1 \ot \psi_2$ with $\psi_1:\ell_1\to \ell$ and $\psi_2:\ell_2\to \ell'$. Then we obtain $\psi = (\psi_1 \ot \psi_2).\id_L$. We deduce that $\psi\sim\id_L$ in $(\Delta_{\mathcal{P}^{op}} U \ot \Delta_{\mathcal{P}^{op}} U')(L)$. 
\epf

\bp \label{tensor-product} (Replacement for \cite[Proposition~5.18]{Moore1})
Let $D$ and $E$ be two $\mathcal{P}$-spaces. Then there is a natural homeomorphism 
\[
\liminj (D \ot E) \iso \liminj D \p \liminj E.
\]
\ep

\bpf
Let $Z$ be a topological space. There is the sequence of natural homeomorphisms 
\[\begin{aligned}
\ttop\big(\liminj (D \ot E),Z\big) & \iso \topdgr\big(D \ot E,\Delta_{\mathcal{P}^{op}}Z\big)\\
&\iso \topdgr \bigg(D,\ell \mapsto \topdgr(E,(s^L_\ell)^*\Delta_{\mathcal{P}^{op}}(Z))\bigg)\\
&\iso \topdgr \bigg(D,\Delta_{\mathcal{P}^{op}}\big(\topdgr(E,\Delta_{\mathcal{P}^{op}}(Z))\big)\bigg)\\
&\iso \ttop\big(\liminj D,\topdgr(E,\Delta_{\mathcal{P}^{op}}(Z))\big)\\
&\iso \ttop\big(\liminj D,\ttop(\liminj E,Z)\big)\\
&\iso \ttop\big((\liminj D) \p (\liminj E),Z \big).
\end{aligned}\]
The proof is complete thanks to the Yoneda lemma.
\epf

Note that in \cite[Theorem~4.3]{Moore2}, the words ``closed symmetric semimonoidal category'' must be replaced by ``biclosed semimonoidal category''.

\section{The case of $\mathcal{G}$-spaces}
\label{non-natural}

\begin{nota}
	In this section, the notations $\ell,\ell',\ell_i,L,\dots$ mean a strictly positive real number.
\end{nota}

For the convenience of the reader, the definition of the reparametrization category $\mathcal{G}$ is recalled: 

\bd Let $\phi_i:[0,\ell_i] \to [0,\ell'_i]$ for $i=1,2$ be two continuous maps preserving the extrema where a notation like $[0,\ell]$ means a segment of the real line. Then the map
\[
\phi_1 \ot \phi_2 : [0,\ell_1+\ell_2] \to [0,\ell'_1 + \ell'_2]
\]
denotes the continuous map defined by 
\[
(\phi_1 \ot \phi_2)(t) = \begin{cases}
\phi_1(t) & \hbox{if } 0\leq t\leq \ell_1\\
\phi_2(t-\ell_1)+\ell'_1 & \hbox{if } \ell_1\leq t\leq \ell_1+\ell_2\\
\end{cases}
\] 
\ed

\begin{nota}
	The notation $[0,\ell_1]\iso^+ [0,\ell_2]$ means a nondecreasing homeomorphism from $[0,\ell_1]$ to $[0,\ell_2]$. It takes $0$ to $0$ and $\ell_1$ to $\ell_2$. 
\end{nota}

\bp \label{paramG} \cite[Proposition~4.9]{Moore1}
There exists a reparametrization category, denoted by ${\mathcal{G}}$, such that the semigroup of objects is the open interval $]0,+\infty[$ equipped with the addition and such that for every $\ell_1,\ell_2>0$, there is the equality \[\mathcal{G}(\ell_1,\ell_2)=\{[0,\ell_1]\iso^+ [0,\ell_2]\}\] where the topology is the compact-open topology (which is $\Delta$-generated by \cite[Proposition~2.5]{Moore2}) and such that for every $\ell_1,\ell_2,\ell_3>0$, the composition map \[\mathcal{G}(\ell_1,\ell_2)\p \mathcal{G}(\ell_2,\ell_3) \to \mathcal{G}(\ell_1,\ell_3)\] is induced by the composition of continuous maps.
\ep

\begin{nota}
	Let $\ell>0$. Let $\mu_{\ell}:[0,\ell]\to [0,1]$ be the homeomorphism defined by $\mu_\ell(t) = t/\ell$. We have $\mu_{\ell}\in \mathcal{G}(\ell,1)$. 
\end{nota}

Recall again that this reparametrization category is not symmetric as a semimonoidal category because the functors $(\ell,\ell')\mapsto \ell+\ell'$ and $(\ell,\ell') \mapsto \ell'+\ell$ coincide on objects, but not on morphisms

\bp \label{dec} Fix $L,\ell_1,\ell_2$. The mapping $(\phi_1,\phi_2) \mapsto \phi_1\ot \phi_2$ induces a continuous bijection which is not a homeomorphism
\[
\displaystyle\bigsqcup\limits_{\substack{\ell'_1>0,\ell'_2>0\\\ell'_1+\ell'_2=L}} \mathcal{G}(\ell'_1,\ell_1)\p \mathcal{G}(\ell'_2,\ell_2) \longrightarrow \mathcal{G}(L,\ell_1+\ell_2)
\]
\ep

\bpf The mapping is a bijection by \cite[Proposition~3.2]{Moore2}. It is continuous since $\mathcal{G}$ is an enriched semimonoidal category. It is not a homeomorphism since the right-hand space is contractible whereas the left-hand one is not.
\epf

\bp \label{B2-continuous}
Fix $L'$. The set map \[B_2:\mathcal{G}([0,2],[0,L'])\longrightarrow \mathcal{G}([0,2],[0,L'])\] which takes $\phi=\phi_1\ot \phi_2$ to $\phi_2\ot \phi_1$ where $\phi_i\in \mathcal{G}([0,1],[0,L'_i])$ with $L'_1=\phi(1)$ and $L'_2=L'-L'_1$ is a idempotent homeomorphism.
\ep

\bpf
It is bijective since $B_2B_2=\id_{\mathcal{G}([0,2],[0,L'])}$. It remains to prove that $B_2$ is continuous. It suffices to prove that $B_2$ is sequentially continuous since the space $\mathcal{G}([0,2],[0,L'])$ is sequential, being metrizable. Let $(\phi^n)_{n\geq 0} = (\phi_1^n\ot\phi_2^n)_{n\geq 0}$ be a sequence of $\mathcal{G}([0,2],[0,L'])$ which converges to $\phi=\phi_1\ot \phi_2$. Then  the sequence $(\phi_i^n)_{n\geq 0}$ converges pointwise to $\phi_i$ for $i=1,2$. Therefore, the sequence $(B_2(\phi^n))_{n\geq 0}$ converges pointwise to $B_2(\phi)$. The proof is complete thanks to \cite[Proposition~2.5]{Moore2}.
\epf

\bp \label{B-continuous}
Fix $L,\ell_1,\ell_2$. There is a unique set map \[B_L^{\ell_1,\ell_2}:\mathcal{G}(L,\ell_1+\ell_2) \to \mathcal{G}(L,\ell_1+\ell_2)\] such that the following diagram of spaces is commutative: 
\[
\xymatrix@C=6em@R=3em
{
	\displaystyle\bigsqcup\limits_{\substack{\ell'_1>0,\ell'_2>0\\\ell'_1+\ell'_2=L}} \mathcal{G}(\ell'_1,\ell_1)\p \mathcal{G}(\ell'_2,\ell_2) \fr{(\phi_1,\phi_2) \mapsto (\phi_2,\phi_1)} \ar@{->}[d]^-{(\phi_1,\phi_2)\mapsto \phi_1\ot \phi_2} & \displaystyle\bigsqcup\limits_{\substack{\ell'_1>0,\ell'_2>0\\\ell'_1+\ell'_2=L}}  \mathcal{G}(\ell'_2,\ell_2) \p \mathcal{G}(\ell'_1,\ell_1)   \ar@{->}[d]^-{(\phi_2,\phi_1)\mapsto \phi_2\ot \phi_1}
	\\
	\mathcal{G}(L,\ell_1+\ell_2) \fr{B_L^{\ell_1,\ell_2}} & \mathcal{G}(L,\ell_1+\ell_2)
}
\]
Moreover, the set map $B_L^{\ell_1,\ell_2}$ is a homeomorphism.
\ep

\bpf
The existence and the uniqueness of the set map $B_L^{\ell_1,\ell_2}$ is a consequence of Proposition~\ref{dec}. It is bijective because all other arrows are bijective. Since $B_L^{\ell_2,\ell_1}B_L^{\ell_1,\ell_2}=\id_{\mathcal{G}(L,\ell_1+\ell_2)}$, it remains to prove that $B_L^{\ell_1,\ell_2}$ is continuous. Observe that 
\[
B_L^{\ell_1,\ell_2}(\psi) = B_2\bigg(\psi \big(\mu^{-1}_{\psi^{-1}(\ell_1)} \ot \mu^{-1}_{L-\psi^{-1}(\ell_1)}\big)\bigg) \big(\mu_{L-\psi^{-1}(\ell_1)} \ot \mu_{\psi^{-1}(\ell_1)}\big).
\]
By \cite[Lemma~6.2]{Moore2}, the mapping $\psi\mapsto \psi^{-1}\mapsto \psi^{-1}(\ell_1)$ is continuous. Thus, the continuity of $B_L^{\ell_1,\ell_2}$ is a consequence of the continuity of $\ot$ proved in Proposition~\ref{dec} and of the continuity of $B_2$ proved in Proposition~\ref{B2-continuous}.
\epf

Let $D$ and $E$ be two $\mathcal{G}$-spaces. The $\mathcal{G}$-space $D\ot E$ is the quotient of 
\[
\bigsqcup_{(\ell_1,\ell_2)} \mathcal{G}(-,\ell_1+\ell_2) \p D(\ell_1) \p E(\ell_2)
\]
by the identifications $(\psi,x_1\phi_1,x_2\phi_2) \sim ((\phi_1 \ot \phi_2)\psi,x_1,x_2)$ by \cite[Corollary~5.13]{Moore1}. Consider the set map 
\[
\bigsqcup_{(\ell_1,\ell_2)} \mathcal{G}(L,\ell_1+\ell_2) \p D(\ell_1) \p E(\ell_2) \longrightarrow (E\ot D)(L)
\]
defined by taking \[(\psi,x_1,x_2)\in \mathcal{G}(L,\ell_1+\ell_2) \p D(\ell_1) \p E(\ell_2)\] to the equivalence class of \[(B_L^{\ell_1,\ell_2}(\psi),x_2,x_1)=(\psi_2\ot \psi_1,x_2,x_1)\] where $\psi=\psi_1\ot \psi_2$ is the unique decomposition of $\psi$ such that $\psi_i\in \mathcal{G}(\ell'_i,\ell_i)$ with $\ell'_1+\ell'_2=L$. It is continuous by Proposition~\ref{B-continuous}. 

The triple $(\psi,x_1\phi_1,x_2\phi_2)$ is taken to the equivalence class of $(\psi_2\ot \psi_1,x_2\phi_2,x_1\phi_1)$. One has $(\phi_1 \ot \phi_2)\psi = (\phi_1\psi_1)\ot (\phi_2\psi_2)$. Therefore, the triple $((\phi_1 \ot \phi_2)\psi,x_1,x_2)$ is taken to the equivalence class of $((\phi_2 \ot \phi_1)(\psi_2\ot \psi_1),x_2,x_1)\sim (\psi_2\ot \psi_1,x_2\phi_2,x_1\phi_1)$. Consequently, we obtain the 

\bth \label{non-natural-homeo}
This mapping yields a continuous map \[B: (D\ot E)(L) \longrightarrow (E\ot D)(L)\] for all $L>0$ and all $\mathcal{G}$-spaces $D$ and $E$ which is a homeomorphism. It is not natural with respect to $L>0$.
\eth

\bpf The map $(D\ot E)(L) \to (E\ot D)(L)$ is not natural with respect to $L\in \Obj(\mathcal{G})$. Indeed, take $(\psi,x_1,x_2)\in (D\ot E)(L)$. Then $(\psi,x_1,x_2)\sim (\id_L,x_1\psi_1,x_2\psi_2)$ in $(D\ot E)(L)$ with $\psi=\psi_1\ot \psi_2$. Consider $\omega:L'\to L$ a map of $\mathcal{G}$. Then $(\id_L,x_1\psi_1,x_2\psi_2)\in (D\ot E)(L)$ is taken to $(\omega,x_1\psi_1,x_2\psi_2) \sim (\id_{L'},x_1\psi_1\omega_1,x_2\psi_2\omega_2)\in (D\ot E)(L')$ with $\omega=\omega_1\ot \omega_2$. On the other hand, $(\id_L,x_2\psi_2,x_1\psi_1)\in (E\ot D)(L)$ is taken to $(\omega,x_2\psi_2,x_1\psi_1)\in (E\ot D)(L')$, and not to $(\omega_2\ot \omega_1,x_2\psi_2,x_1\psi_1)$.  
\epf

Note that the mapping $(\psi,x_1,x_2)\mapsto (\psi,x_2,x_1)$ does not induce a map from $(D\ot E)(L)$ to $(E\ot D)(L)$. Indeed, $(\psi,x_1\phi_1,x_2\phi_2)$ is taken to $(\psi,x_2\phi_2,x_1\phi_1)$ whereas $((\phi_1\ot \phi_2)\psi,x_1,x_2)$ is taken to $((\phi_1\ot \phi_2)\psi,x_2,x_1)$ and $(\psi,x_2\phi_2,x_1\phi_1) \sim ((\phi_2\ot \phi_1)\psi,x_2,x_1)$ which is not equal to $((\phi_1\ot \phi_2)\psi,x_2,x_1)$ in $(E\ot D)(L)$ in general.

%\bibliographystyle{../plainurlwithoutprefixDOI} 
%%\bibliographystyle{alpha} 
%{\footnotesize
%	\bibliography{../Bibliotheque}
%}

\end{document}